\documentclass[a4paper]{article}
\usepackage{amsmath, graphicx, amssymb}
\usepackage[utf8]{inputenc}
\newcommand{\pb}{\,\begin{matrix} \\ +\end{matrix}\,}
\title{Ramanujan in Computing Technology\footnote{A preliminary version of the paper was presented in the National Conference on Advances in Computing Technology 2020 (NC-ACT 2020) organised by Department of Computer Applications, Vidya Academy of Science and Technology, Thrissur - 680501, Kerala, India during 4 - 5 December 2020.}\\[5mm]}
\author{{\bf V N Krishnachandran} \\ Vidya Academy of Science and Technology\\ (APJ Abdul Kalam Technological University)\\ Thrissur - 680501, Kerala, India \\
(Email: {\tt krishnachandranvn@gmail.com})}
\date{}
\begin{document}
\maketitle
\tableofcontents
\newpage
\begin{center}
\includegraphics[height=5cm]{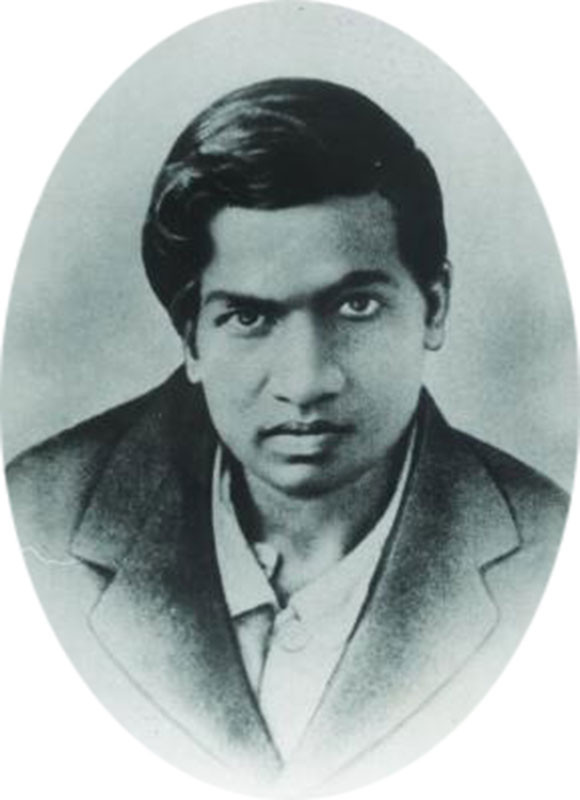}

Srinivasa Ramanujan

(22 December 1887 - 26 April 1920)
\end{center}
\section{Introduction}
This paper is a tribute  to the genius of the legendary Indian mathematician Srinivasa Ramanujan (22 December 1887 - 26 April 1920) in the centenary year of his death. The life story of Ramanujan is so well known that it needs no elaboration not even a summarisation. In his short life period he made substantial contributions to mathematical analysis, number theory, infinite series, and continued fractions, including solutions to mathematical problems then considered unsolvable. Ramanujan independently compiled nearly 3,900 results in the form of identities and equations. Many were completely novel; his original and highly unconventional results, such as the Ramanujan prime, the Ramanujan theta function, partition formulae and mock theta functions, have opened entire new areas of work and inspired a vast amount of further research. Nearly all his claims have now been proven correct.

The focus of the paper is the increasing influence of the ideas propounded by Ramanujan in the development of computing technology. We shall discuss the application of certain infinite series discovered by Ramanujan in computing the value of the mathematical constant $\pi$. We shall also consider certain special graphs known as Ramanujan graphs and the reason for designating them as such. We shall examine how certain researchers are attempting to create an abstract machine which they call Ramanujan machine which is thought of as simulating the  hypothesised thought process of Ramnujan. We shall also have a brief look at the applications of Ramanujan's discoveries in signal processing. 
\section{Some milestones in the history of computation of $\pi$}
The value of the ratio of the circumference of a circle to its diameter has been of interest to mankind since the beginning of civilisations. However, the modern notation for the value of the ratio of the circumference of a circle to its diameter, namely $\pi$, was first used in print in a work of an English mathematician William Jones (1675-1749) published in 1706. However, it became a widely accepted notation only when Euler used it in his famous 1748 publication "{\em Introductio in analysin infinitorum}" (Introduction to the Analysis of the
Infinite). 

\subsection{Ancient times}

The Rhind Papyrus dated around 1650 BCE says: 
``Cut off 1/9 of a diameter and construct a square upon the remainder; this has the same area as the circle". This gives a value of 
$$
4(8/9)^2 = 3.16049
$$
to the ratio of the circumference of a circle to its diameter, which is fairly accurate. The {\em Shulba Sutras} composed during the period 800 - 200 BCE  and considered to be appendices to the Vedas  give several different approximations to the ratio the circumference of a circle to its diameter and among them the value closest to the actual value is 3.088 (see \cite{HistoryOfPi}). In the earliest existing Chinese mathematical texts, {\em Zhoubi suanjing} (The Mathematical Classic of the Zhou Gnomon), a book dated to the period of the Zhou dynasty (1046 - 256 BCE) and {\em Jiuzhang suanshu} (Nine Chapters on the Mathematical Art), the ratio of the circumference of the circle to its diameter was taken to be three (see \cite{China}). 

\subsection{Geometric algorithms}
\begin{figure}[!h]
\begin{center}
\includegraphics[width=6cm]{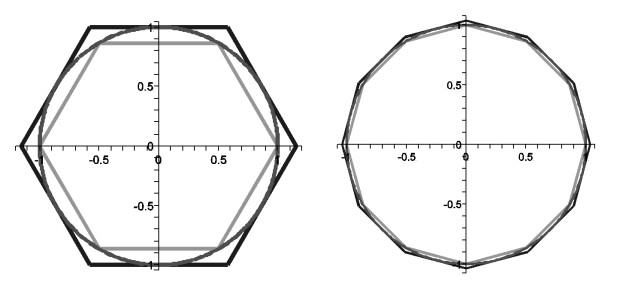}
\caption{Archimedes' idea for computing the value of $\pi$}
\end{center}
\end{figure}
The first algorithm for calculating the value of the ratio of the circumference of a circle to its diameter was a geometrical approach using polygons, devised around 250 BCE by the Greek mathematician Archimedes. Archimedes computed upper and lower bounds of the value of this ratio by drawing a regular hexagon inside and outside a circle, and successively doubling the number of sides until he reached a 96-sided regular polygon. By calculating the perimeters of these polygons, he proved that the ratio lies between $223/71$ and $22/7$. Mathematicians using polygonal algorithms calculated $39$ digits of the value of the ratio in 1630, a record only broken in 1699 when infinite series were used to reach 71 digits (see \cite{Unleashed}).

\begin{figure}[!h]
\begin{center}
\includegraphics[height=7cm]{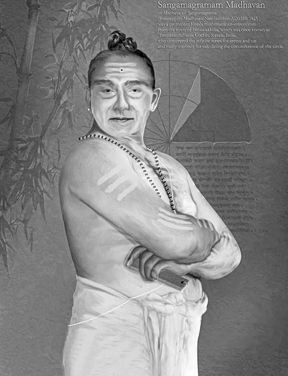}
\caption{Sangamagrama Madhava (c.1340 – c.1425). This is a digital image of Madhava drawn up with inputs provided by his descendants and released in 2014 by the Madhava Ganitha Kendram, a Kochi based voluntary association working to revive his works.}\label{Mad}
\end{center}
\end{figure}
 
\subsection{Sangamagrama Madhava}
Sangamagrama Madhava (c.1340 - c.1425), an Indian mathematician and astronomer considered the founder of the Kerala school of astronomy and mathematics, calculated the value of $\pi$ correct to $11$ decimal places as 
$$
\pi= 3.14159265359.
$$ 
Madhava obtained this value probably by taking the  first 21 terms in the following infinite series discovered by him (see \cite{Gupta}):
$$
\pi 
={\sqrt {12}}\sum _{k=0}^{\infty }{\frac {(-3)^{-k}}{2k+1}}
={\sqrt {12}}\left(1-{1 \over 3\cdot 3}+{1 \over 5\cdot 3^{2}}-{1 \over 7\cdot 3^{3}}+\cdots \right)
$$
\subsection{William Shanks}
The British amateur mathematician William Shanks spending over 20 years attempted to calculate $\pi$ to $707$ decimal places. He did calculate 707 digits, but only the first $527$ were correct. This was accomplished in 1873 and  this was the longest expansion of $\pi$ until the advent of the electronic digital computer.

\begin{figure}[!h]
\begin{center}
\includegraphics[height=4cm]{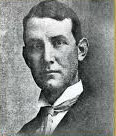}
\caption{William Shanks (1812 - 1882)}\label{Shanks}
\end{center}
\end{figure}

To compute the value of $\pi$ Shanks used the following result known as Machin formula (discovered by  John Machin, an English astronomer, in 1706)
$$
\frac{\pi}{4}=4 \tan^{-1}\left(\frac{1}{5}\right)-\tan^{-1}\left(\frac{1}{239}\right),
$$
and the Maclaurin's series expansion 
$$
\tan^{-1}x=x - \frac{x^3}{3}+\frac{x^5}{5}-\frac{x^7}{7}+\cdots
$$
This formula has a significantly increased rate of convergence, which makes it a much more practical method of calculation and it remained as the primary tool of $\pi$ calculations for centuries (well into the computer era). Machin himself used this formula to compute $\pi$ to $100$ decimal places.

\subsection{Computer programmes}
 Shank's computation is only child's play in the computer era as illustrated by the fact that the following code written in the C language prints out as much as $800$ digits of $\pi$ accurately (for a detailed analysis of why the program does as claimed, see \cite{pi800}):
\begin{verbatim}
   int a=10000,b,c=2800,d,e,f[2801],g;
   main()
   {
     for( ;b-c;) f[b++]=a/5;
     for( ;d=0,g=c*2;c=14,printf("%.4d",e+d/a),e=d%a)
        for(b=c; d+=f[b]*a,f[b]=d%--g,d/=g--,--b;d*=b);
   }
\end{verbatim}
It is even simpler if we use a computer algebra system  to get 800 digits of $\pi$. For example, in Maxima we need only issue the following commands to get 800 digits of $\pi$:
\begin{verbatim}
   fpprintprec:800$
   set_display(ascii)$
   bfloat(%pi);
\end{verbatim}
The current world record for the number of calculated digits of $\pi$ is $50$ trillion and it was created by Timothy Mullican on 29 January 2020 (see \cite{ycruncher}). 
\subsection{Do we ever need all these digits?}
The largest number of digits of $\pi$ that we will ever need is 42, at least for computing circumferences of circles (see \cite{piday}). To compute the circumference of the known universe with an error less than the diameter of a proton, we need only 42 digits of $\pi$,  assuming that the diameter of the known universe is $93$ billion light years and that the diameter of a proton is $1.6\time 10^{-15}$ metres.  Thus, in the fifty trillion digits of $\pi$ computed for the current record, all digits beyond the 42nd have no practical value.
Here are all the digits of $\pi$ we will ever need:
$$
\pi = 3.141592653589793238462643383279502884197169.
$$

\section{Ramanujan's series for computing $\pi$}
In the long journey covering the milestones in the computation of digits of $\pi$, Srinivasa Ramanujan's name appears sometime around the late 1970's and early 1980's and his name appears via a much celebrated infinite series known as the ``Ramanujan series''. 
\subsection{Ramanujan's series}
\begin{figure*}[!h]
\begin{center}
\includegraphics[width=10cm]{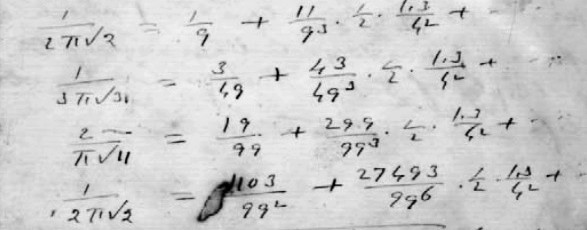}
\caption{The series used by Bill Gosper as it appears in Ramanujan's Notebook 3}\label{Series1}
\end{center}
\end{figure*}

\begin{figure*}[!h]
\begin{center}
\includegraphics[width=12cm]{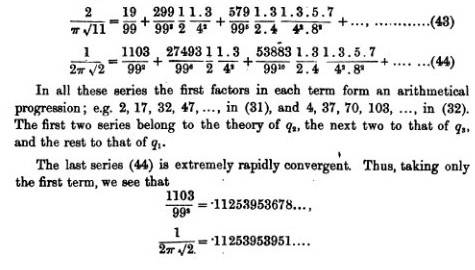}
\caption{The series used by Bill Gosper (see equation (44)) and comments on it in Ramanujan's paper (see \cite{Ramanujan1914})}\label{Series2}
\end{center}
\end{figure*}

In 1903, or perhaps earlier, while in school, Ramanujan began 
to record his mathematical discoveries in notebooks. He continued this practice till his departure to England in 1914. There are three such notebooks the originals of which are now preserved in Madras University, Chennai. A fourth notebook, which had been thought to have been lost and hence referred to as ``The Lost Notebook'' was later unearthed in 1974. It turned out to be not a notebook but a collection of loose sheets of paper. In the last page of the third of these notebooks Ramanujan listed a large number of series expansions for $\frac{1}{\pi}$ (see Figure \ref{Series1}). 
After his arrival in England, Ramanujan published a 
paper in 1914 (see \cite{Ramanujan1914}) containing all these  expansions for $\frac{1}{\pi}$ (see Figure \ref{Series2}). Interestingly all of them remained unproven for nearly half a century. 

The following is the series given by Ramanujan and used by Gosper to compute the digits of $\pi$:
$$
{\displaystyle {\frac {1}{\pi }}={\frac {2{\sqrt {2}}}{99^{2}}}\sum _{k=0}^{\infty }{\frac {(4k)!}{k!^{4}}}{\frac {26390k+1103}{396^{4k}}}} 
$$

\subsection{Using Ramanujan's series for $\pi$ computation}
In 1977, Bill Gosper took the last of Ramanujan’s series from the list referred to, and used it to compute a record number of digits of $\pi$.  There soon followed other computations, all based directly on Ramanujan’s idea. 

Ralph William Gosper Jr (born April 26, 1943), known as Bill Gosper, is an American mathematician and programmer. Gosper was one among the first few persons to realize the possibilities of symbolic computations on a computer as a mathematics research tool whereas computer methods were previously limited to purely numerical methods. He made major contributions to Macsyma, Project MAC's computer algebra system. 
\begin{figure}
\begin{center}
\includegraphics[width=4cm]{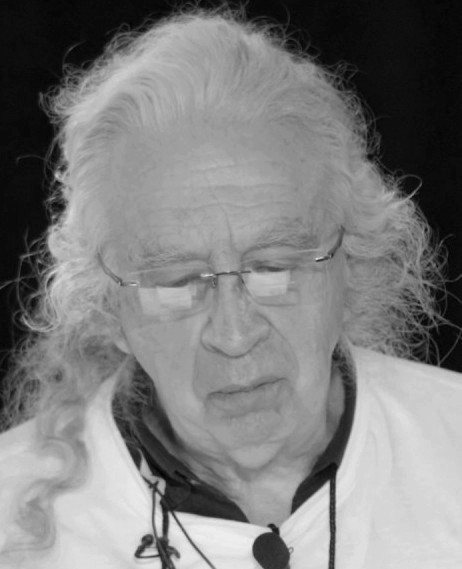}
\caption{Bill Gosper}
\end{center}
\end{figure}

In  November, 1985, Gosper employed a lisp machine at Symbolics 
and Ramanujan's 
series to calculate 17,526,100 digits of $\pi$, which at that 
time was a world record (see \cite{maa}). It was only 1n 1987, the Ramanujan series for $1/\pi$ was finally 
proved and the proof was accomplished by Jonathan and Peter 
Borwein (see \cite{maa}). 

\subsection{Comments on the computations}

There  were  a  few  interesting  things about Gosper's
 computation.To put these computations in perpective, it is better to quote verbatim from a reference document on $\pi$ calculations published by Princeton University (see \cite{HistoryofPiCalculations}):

\begin{itemize}
\item
First, when he decided to use  that particular formula, there was no
proof that it actually converged to $\pi$!  Ramanujan  never  gave  the
math  behind  his  work,  and  the Borweins had not yet been able to
prove it, because there was some  very  heavy math that needed to be
worked through.  It  appears  that  Ramanujan  simply  observed  the
equations  were  converging  to  the  1103  in the formula, and then
assumed it must actually be 1103.  (Ramanujan was not known  for
rigor  in his math, or for providing any proofs or intermediate math
in his formulas.) The  math  of  the  Borwein's  proof was such that
after he had computed 10 million digits, and verified them against a
known  calculation,  his  computation  became  part  of  the  proof.
Basically it was like, if you have two integers  differing  by  less
than one, then they have to be the same integer.
\item
The second interesting  thing  is  that  he  chose  to  use continued
fractions to do his  calculations.   Most  calculations  before  and
since  were  done  by  direct  calculation to the desired precision.
Before you did any calculations,  you  had to decide how many digits
you wanted, and later if you wanted more, you had to start over from
the beginning.  By using continued fractions,  and  a  novel  coding
style,  he was able to make his resumable.  He could add more digits
any time he felt like it and had the spare time.  This was  a  major
breakthrough  at  the  time,  because  all previous efforts required
starting over from the beginning if you wanted more.
\item
The  third  interesting  thing  about his calculations was the other
reason he chose to use infinite simple continued fractions.   It  is
still  not  known  whether pi has any 'structure' or patterns to it.
It is known that  it's  irrational  and  transcendental, but it still
might have some pattern to  it  that  would  allow  us  more  easily
calculate  its  digits.   We  just don't know.  And patterns show up
more readily as a continued  fraction  rather than in some arbitrary
base that we humans call 'base 10'.  As  an  example,  'e'  and  the
square  root of two both have very simple, obvious patterns in their
continued fractions, even though they  appear to be non-repeating in
base 10.
\end{itemize}

\section{Chudnovsky algorithm}
\subsection{Chudnovski's formula for computing $\pi$}
Further developing the ideas on which the Ramanujan series was proved, Chudnovsky brothers (David Volfovich Chudnovsky (born 1947) and Gregory Volfovich Chudnovsky (born 1952)) developed another series for $1/\pi$: 
$$
 {\displaystyle {\frac {1}{\pi }}=12\sum _{q=0}^{\infty }{\frac {(-1)^{q}(6q)!(545140134q+13591409)}{(3q)!(q!)^{3}\left(640320\right)^{3q+3/2}}}}
$$

This formula was used in the world record calculations of 2.7 trillion digits of $\pi$ in December 2009, 10 trillion digits in October 2011, 22.4 trillion digits in November 2016, 31.4 trillion digits in September 2018 - January 2019 (see \cite{Google}), and 50 trillion digits on 29 January  2020 (see \cite{ycruncher}). 

The computer algebra system {\em Matheamtica} widely used in many scientific, engineering, mathematical, and computing fields uses this formula to compute $\pi$ (see \cite{Unleashed}, \cite{Wolfram}).

\subsection{Chudnosky brothers}
\begin{figure}[!h]
\begin{center}
\includegraphics[width=6cm]{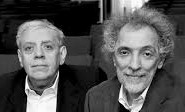}
\caption{Chudnovski brothers}
\end{center}
\end{figure}
David Volfovich Chudnovsky and Gregory Volfovich Chudnovsky are American mathematicians and engineers known for their world-record mathematical calculations and developing the Chudnovsky algorithm used to calculate the digits of $\pi$ with extreme precision. They were born in Kiev in Ukraine which at the time of their births was part of the erstwhile Soviet Union. In part to avoid religious persecution and in part to seek better medical care for Gregory, who had been diagnosed with myasthenia gravis, a neuromuscular disease, the Chudnovsky family emigrated to the United states in 1977-78 and settled in New York. Gregory works on problems of diophantine geometry and transcendence theory, high-performance computing and computing architecture, and mathematical physics and its applications. A 1992 article in The New Yorker quoted the opinion of several mathematicians that Gregory Chudnovsky was one of the world's best living mathematicians (see \cite{NewYorker}).
\subsection{The computational procedure}
For a high performance iterative implementation, Chudnovski's series may have to be put in different forms. We present below one of these approaches (see \cite{Wiki}). In this approach, firstly,  the series is recast in the following form.
$$
{\frac {426880{\sqrt {10005}}}{\pi }}=\sum _{q=0}^{\infty }{\frac {(6q)!(545140134q+13591409)}{(3q)!(q!)^{3}\left(-262537412640768000\right)^{q}}}
$$

There are three big integer terms in the series, namely, a multinomial term $M_q$, a linear term $L_q$, and a exponential term $X_q$ as shown below:
$$
 {\displaystyle \pi =C\left(\sum _{q=0}^{\infty }{\frac {M_{q}\cdot L_{q}}{X_{q}}}\right)^{-1}},
$$
 where
\begin{align*}
\displaystyle C &=426880{\sqrt {10005}}\\
\displaystyle M_{q}&={\frac {(6q)!}{(3q)!(q!)^{3}}}\\
\displaystyle L_{q}& =545140134q+13591409\\
\displaystyle X_{q}&=(-262537412640768000)^{q}.
\end{align*}
The terms $M_q$, $L_q$, and $X_q$ can be computed by using the following recurrence relations
\begin{align*}
L_{q+1}&=L_{q}+545140134\\
X_{q+1}&=X_{q}\cdot (-262537412640768000)\\
M_{q+1}&=M_{q}\cdot \left({\frac {(12q+2)(12q+6)(12q+10)}{(q+1)^{3}}}\right)
\end{align*}
with
$$
L_{0}=13591409, \quad X_0=1, \quad M_0=1.
$$
The computation of $M_q$ can be further optimized by introducing an additional term $K_q$ as follows:
\begin{align*}
K_{q+1}&=K_{q}+12\\
M_{q+1}&=M_{q}\cdot \left({\frac {K_{q}^{3}-16K_{q}}{(q+1)^{3}}}\right)
\end{align*}
with
$$
K_0=6, \quad M_0=1.
$$
\section{Ramanujan graphs}
Ramanujan graphs are graphs having a very special property. 
The terminology ``Ramanujan graph''  was first coined by Lubotzky, Phillips and Sarnak in a paper published in 1988 (see \cite{RGraph}). Though Ramanujan had no special interest in graph theory and these graphs are not due to Ramanujan, they are named after Srinivasa Ramanujan because the  Ramanujan-Petersson conjecture is used in proving that an important  class of graphs have the special property. First let us see what are Ramanujan graphs.

Let $G=(V,E)$ be an undirected graph, $V$ being the set of vertices and $E$ the set of edges. The degree of a vertex $v$ in $G$ is the number of edges connecting it to other vertices in the graph. We say that a graph $G$ is $k$-regular  if every vertex $v  G$ has degree $k$. A path in a graph is a sequence of edges which joins a sequence of vertices. A graph is connected if  there is a path between every pair of vertices. 

We can associate a matrix called an adjacency matrix with every graph $G$. Let $V=\{v_1, \ldots, v_n\}$. The adjacency matrix of $G$ is the $n\times n$ matrix $A=[a_{ij}]$ where
$$
a_{ij}=
\begin{cases}
1 & \text{ if an edge joining $v_i$ and $v_j$ exists}\\
0 & \text{ otherwise.}
\end{cases}
$$
Since $G$ is undirected, $A$ is a symmetric matrix and all of its eigenvalues are real. Moreover, if $G$ is $k$-regular, the largest eigen value of $A$ is $k$. 
\subsection{Definition}
A $k$-regular connected graph $G$  is called a Ramanujan graph\footnote{The terminology {\em Ramanujan Graph} was suggested by Peter Sarnak (see \cite{Sarnak}).} if
the second largest (in absolute value) eigenvalue of the adjacency matrix $A$ of $G$, denoted by $\lambda(G)$, satisfies the property
$$\lambda(G) \le 2\sqrt{k-1}.$$

\begin{figure}
\begin{center}
\includegraphics[width=2.5cm, height=2.5cm]{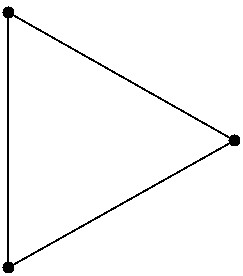}
\qquad
\includegraphics[width=2.5cm, height=2.5cm]{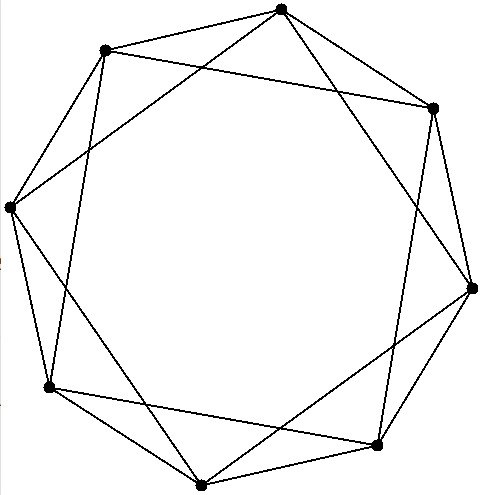}
\caption{Examples of Ramanujan Graphs}\label{RG}
\end{center}
\end{figure}

The graphs shown in Figure \ref{RG} are examples of Ramnujan graphs.

\subsection{Ramanujan conjecture}
Ramanujan studied the function $\Delta(z)$ defined by
\begin{align*}
\Delta(z)%
&=q\prod_{n>0}\left (1-q^n \right)^{24}\\
&= q-24q^2+252q^3- 1472q^4 + 4830q^5-\cdots\\
&= \sum_{n>0}\tau(n)q^n
\end{align*}
where $ q=e^{2\pi iz}$.  The function $\tau(n)$ is called the Ramanujan $\tau$-function. Ramanujan computed the values of $\tau(p)$ when $p$ is assigned the first $30$ prime numbers as values. These values are shown in Table \ref{tauvalues}. 
\begin{table*}[!h]
\begin{center}
\includegraphics[width=12cm]{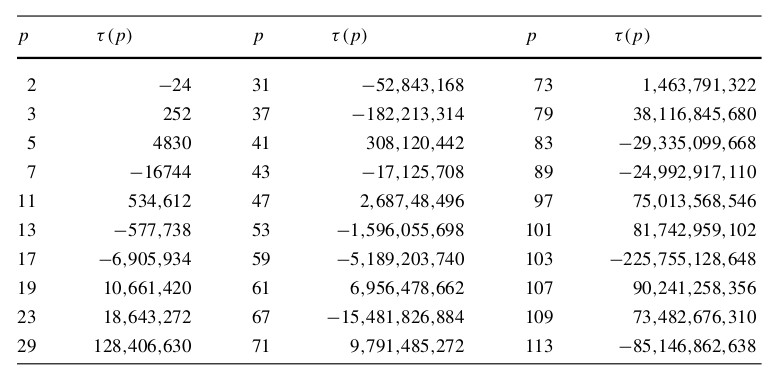}
\caption{Values of $\tau(p)$ calculated by Ramanujan}\label{tauvalues}
\end{center}
\end{table*}

Based on these values of $\tau(p)$ Ramanujan made the following  observations regarding the properties of the function $\tau(n)$ (see \cite{tau}):
\begin{enumerate}
\item
If $m$ and $n$ are relatively prime, then $$\tau(mn)=\tau(m)\tau(n).$$ 
\item
For a prime number $p$ and an integer $j$ we have: 
$$
\tau(p^{j+1})=\tau(p)\tau(p^j)-p^{11}\tau(p^{j-1}).
$$
\item
For  a prime number $p$ we have:
$$
|\tau(p)| \leq 2p^{\frac{11}{2}}.
$$
\end{enumerate}

 Ramanujan did not prove these results. In 1917, L. Mordell proved the first two observations using techniques from complex analysis. The third observation remained unproven for a long time and it came to be referred to as the {\em Ramanujan conjecture} (see \cite{Murti}).
It was finally proven  in 1974 by Deligne\footnote{Pierre  Deligne is a Fields medal winning  Belgian mathematician best known for his work on the Weil conjectures.  For a layman's introduction to Deligne's work, one may read \cite{Gowers}.}.

Let $p, q$ be  prime numbers and $x_1, x_2, x_3, x_4$ be integers. Consider the expression
$$
Q(x_1, x_2, x_3, x_4) = x_1^2 + (2qx_2)^2 +(2qx_3)^2 + (2qx_4)^2
$$
Given a nonnegative integer $n$, let  $r_Q(n)$ denote the number of representations of $n$ in the form $Q(x_1, x_2, x_3, x_4)$. There is no explicit formula for $r_Q(n)$. But using the Ramnujan conjecture we get a good approximation for $r_Q(n)$ in the form
$$
r_Q(p^k) = C(p^k) + O(p^{k(1/2+\epsilon)})
$$
where 
$$
C(p^k) = c
\frac{p^{k+1}-1}{p-1}
$$
for some constant $c$, and $O$ denote the well known Big-O notation.
\section{Construction of Ramanujan graph}
\subsection{Cayley graph}
 Let $G$ be a group and $S$ be some symmetric set of elements
of G (that is, $S$ is a subset of $G$ such that if $x\in S$ then $x^{-1}\in S$). Then the Cayley graph generated by $(G, S)$ is the graph whose vertices are the elements of $G$, and given $g, h \in G$, the edge $(g, h)$ exists if and only if there
exists some $s \in S$ such that $gs = h$.

\begin{figure}
\begin{center}
\begin{tabular}{cc}
\includegraphics[width=4cm]{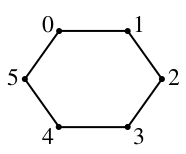}\qquad\qquad&
\includegraphics[width=4cm]{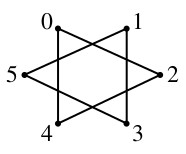}\\
(a) $S=\{-1,1\}$ & (b) $S=\{-2,2\}$\\
\end{tabular}
\caption{Cayley graphs with the group $\mathbb Z_6$}\label{graph3}
\end{center}
\end{figure}

As an example, consider the group whose set of elements is  
$$
\mathbb Z_6=\{0,1,2,3,4,5\}
$$
and the binary operation is addition modulo $6$. If we choose $=\{-1,1\}$ we get the graph in Figure \ref{graph3}(a) and if we set $S=\{-2,2\}$ we get the graph in Figure \ref{graph3}(b).

\subsection{The groups $PGL(2,q)$ and $PSL(2,q)$}
Let $q$ be a prime number. 
We denote by $Z_q$ the field whose set of elements is the set $\{0,1,2, \ldots, q-1\}$, the addition operation is addition modulo $q$ and multiplication operation is multiplication modulo $q$. Let $Z_q^2$ be the 2-dimensional vector space over $Z_q$ and $GL(2,q)$  be the set of nonsingular linear transformations of $Z_q^2$ ($GL$  for {\em General Linear}). Each such linear transformation can be represented by a square matrix of order $2$ whose elements are in $Z_q$ and whose determinant is nonzero (modulo $q$). The set $GL(2,q)$ is a group under operation of multiplication of matrices. 

For $A,B\in GL(2,q)$, we write $A\sim B$ if there is a scalar matrix $C\in GL(2,q)$ such that $A=BC$. (A scalar matrix is a diagonal matrix  with equal-valued elements along the diagonal.) This is an equivalence relation in $GL(2,q)$ and the set of equivalence classes is denoted by $PGL(2,q)$ ($PGL$ for {\em }Projective General Linear). It can be easily seen that $PGL(2,q)$ is also a group under an obvious multiplication operation. 

The subset of $PGL(2,q)$ consisting of matrices with determinant equal to $1$ is also a group and is denoted by $PSL(2,q)$ ($PSL$ for {\em Projective Special Linear}). 
\subsection{A set for generating a Cayley graph of $PGL(2,q)$ and $PSL(2,q)$}
Let us now assume that $p,q$ satisfy the conditions
$$
p\equiv 1\,\, (\text{mod } 4), \quad q \equiv 1 \,\,(\text{mod }4).
$$
The second of these conditions imply that there is an integer $i$ such that
$$
i^2\equiv -1\,\, (\text{mod } q).
$$
We choose such an $i$ and keep it fixed.

Next we consider the equation:
$$
a_0^2+a_1^2+a_2^2+a_3^2 =p.
$$
The condition $p\equiv 1\,\, (\text{mod } 4)$ implies that there are precisely $p+1$ solutions to this equation with $a_0>0$ and $a_1, a_2, a_3$ even. We now form  the set $S$ of $p+1$ matrices with the following form using these $p+1$ solutions:
$$
\alpha =\begin{bmatrix}\phantom{-}a_0+ia_1 & a_2+ia_3\\-a_2+ia_3 & a_0-ia_1\end{bmatrix}
$$
It can be seen that $S$ is a symmetrical set in $PGL(2, q)$. 
\subsection{The graph $X^{p,q}$}
Let $m,n$ be positive integers. We say that {\em $m$ is a quadratic residue $\text{mod } n$} if there is an integer $x$ such that
$$
x^2\equiv m\,\, (\text{mod } n).
$$
If there is no such integer $x$ we say that {\em $m$ is a quadratic nonresidue $\text{mod }n$}.
For example, $1,3,4,5,9$ are quadratic residues mod $11$ and $2,6,7,8,10$ are quadratic nonresidues mod $11$. 
The graph $X^{p,q}$ is defined as follows:

\begin{itemize}
\item
If $p$ is a quadratic residue mod $q$, then  the Cayley
graph of 
$$
PSL(2, q)
$$ 
with $S$ as described above as the generating set is defined as $X^{p,q}$. In this case, we have a $(p + 1)$-regular graph on $q(q^2 - 1)/2$ vertices. 
\item
If $p$ is a quadratic non-residue mod $q$, then  the Cayley
graph of 
$$
PGL(2, q)
$$ 
 with $S$ as described above as the generating set is defined as $X^{p,q}$. This   gives a $(p + 1)$-regular graph on $q(q^2 - 1)$ vertices. 
\end{itemize}
%

\subsection{$X^{p,q}$ is a Ramanujan graph}

A. Lubotzky, R. Phillips and P. Sarnak, in 1988 (see \cite{RGraph}) proved that $X^{p,q}$ is a Ramanujan graph, that is, 
$$
\lambda(X^{p,q})\le 2 \sqrt{p+1}.
$$

The proof of this required an application of the Ramanujan conjecture. So, A. Lubotzky, R. Phillips and P. Sarnak named any graph having such a property as a Ramanujan graph. 
\subsection{Numerical example}

First of all we have to choose two prime numbers satisfying the following conditions:
\begin{itemize}
\item
$p\equiv 1 (\text{mod } 4)$
\item
$q\equiv 1 (\text{mod } 4)$
\item
$p$ is a quadratic residue mod $q$, that is, there exists an integer $x$ such that $x^2\equiv p (\text{mod } q)$. 
\end{itemize}
It can be shown that smallest such prime numbers are $p=5$ and $q=29$. Note that $11^2\equiv 5\,\, (\text{mod } 29)$.

Next we consider the diophantine equation
$$
a_1^2 + a_2^2 + a_3^2 +a_4^2=5.
$$
The $p+1 =6$ solutions to this equation with $a_0>0$ and $a_1, a_2, a_3$ even are
$$
\{(1,\pm 2, 0, 0), (1,0,\pm 2,0), (1,0,0,\pm2)\}.
$$
We have to choose an integer $i$ such that $i^2\equiv -1 \,\,(\text{mod } q)$. We may choose $i=12$ 
because $12^2 \equiv -1 \,\,(\text{mod } 29$. (A second possible value of $i$ is $17$.) 

The elements of the generating set $S$ of $PGL(2,29)$ are given by
$$
\begin{bmatrix}
\phantom{-}a_0+ia_1 & a_2+ia_3\\-a_2+ia_3 & a_0-ia_1
\end{bmatrix}
$$
Substituting for $a_0,a_1,a_2,a_3$ and $i$, and making computations modulo $29$ we get the following generating set for $PGL(2,29)$:
$$
S=\left\{
\begin{bmatrix} 25 & 0 \\ 0 & 6 \end{bmatrix}, 
\begin{bmatrix} 6 & 0 \\ 0 & 25\end{bmatrix}, 
\begin{bmatrix} 1 & 2\\ 27 & 1\end{bmatrix},
\begin{bmatrix} 1 & 27 \\ 2 & 1\end{bmatrix}, 
\begin{bmatrix} 1 & 24\\ 24 & 1\end{bmatrix},
\begin{bmatrix} 1 & 5\\ 5 & 1\end{bmatrix}
\right\}.
$$

To get a generating set for $PSL(2,29)$, we need a set of matrices whose determinants are equal to $1$ modulo $29$. To get such a set we have to multiply the elements of the set $S$ by $\frac{1}{\sqrt{p}}\,\,(\text{mod } q$. There are two values for $\sqrt{5}\,\,(\text{mod } 29$, namely, $11$ and $18$. We choosing $11$ arbitrarily, we have 
$$
\frac{1}{\sqrt{5}}\equiv \frac{1}{\sqrt{11}}\equiv 8\,\, (\text{mod } 29).
$$
Multiplying the elements in $S$ by $8$ modulo $29$ we get the following generating set for $PSL(2,29)$.
$$
S^\prime=\left\{
\begin{bmatrix} 26 & 0 \\ 0 & 19 \end{bmatrix}, 
\begin{bmatrix} 19 & 0 \\ 0 & 26\end{bmatrix}, 
\begin{bmatrix} 8 & 16\\ 13 & 8\end{bmatrix},
\begin{bmatrix} 8 & 13 \\ 16 & 8\end{bmatrix}, 
\begin{bmatrix} 8 & 18\\ 18 & 8\end{bmatrix},
\begin{bmatrix} 8 & 11\\ 11 & 8\end{bmatrix}
\right\}.
$$

The number of elements in the set $PSL(5,29)$ is
$$
\tfrac{1}{2}q(q^2-1)=\tfrac{1}{2}\times 29(29^2-1)=12180.
$$
Thus the Ramanujan graph $X^{5,29}$ is a $6$-regular graph with $12180$ vertices!
\section{Applications of Ramanujan graphs}
Most applications of Ramanujan graphs are a consequence of the fact that Ramanujan graphs are good examples of a certain class of graphs known as {\em expander graphs}. So we give here a quick review of the concept of expander graphs.

\subsection{Expander graphs}
\begin{figure}
\begin{center}
\includegraphics[width=5cm]{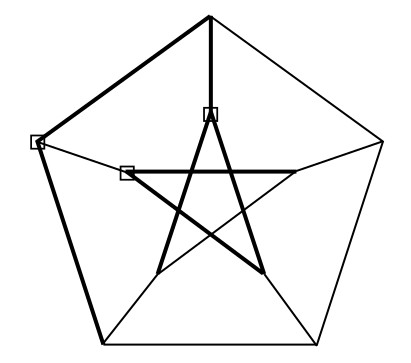}
\caption{A set of vertices (surrounded by squares) and its boundary (thick edges)}\label{Boundary}
\end{center}
\end{figure}

Let $X=(V,E)$ be an undirected connected graph where $V$ is the set of vertices and $E$ is the set of edges. Let $F$ be a set of vertices in $X$. The {\em boundary} $\partial F$ is the set of edges connecting vertices in $F$ to vertices in $V-F$. For example, in the graph shown in Figure \ref{Boundary} let $F$ be the verices surrounded by squares. The set of edges shown by thick line segments is the boundary $\partial F$ of $F$.

The {\em expanding constant}, or {\em isoperimetric constant} of $X$, is
$$
h(X)=\inf\left\{\frac{|\partial F|}{\min \{|F|,|V-F|\}}:
F\subseteq V,\,\, 0<|F|<+\infty\right\}.
$$

If we view $X$ as a network transmitting information (where information retained by some vertex propagates, say in $1$ unit of time, to neighboring vertices), then $h(X )$ measures the ``quality'' of $X$ as a network: if $h(X )$ is large, information propagates well. For example, if $X=K_n$ is the highly connected complete graph on $n$ vertices then $h(X)\sim \frac{n}{2}$ which grows proportionally with the number of edges. On the other hand if $X=C_n$ is the minimally connected cyclic graph on $n$ vertices then $H(X)\le \frac{4}{n}$ which tends to zero as the $n$ increases to infinity. Thus $h(X)$ is a measure of the connectivity of $X$ as a network.

An expander family of graphs is defined as follows.

Let $(X_m )_{m\ge 1}$ be a family of graphs $X_m = (V_m , E_m )$ indexed
by $m \in \mathbb N$ (set of natural numbers). Furthermore, fix $k \ge  2$. Such a family $(X_m )_{ m\ge 1}$ of finite, connected,
$k$-regular graphs is a family of expanders if $|V_m | \rightarrow +\infty$ for $m \rightarrow +\infty$, and
if there exists $\epsilon  > 0$, such that $h(X_m ) \ge \epsilon$ for every $m \ge 1$.

\subsection{Network theory}
Consider the construction
of a new telephone network and 
we want the network to have a high degree of connectivity. However, laying lines will be expensive, so we want
to achieve this high degree of connectivity using as few edges in the network as possible. This problem can be modeled as one related to expander graphs. An expander graph is a sparse graph (a graph having only a few edges)  that has strong connectivity properties, quantified using vertex, edge or spectral expansion. Ramanujan graphs has applications in the construction and study of expander graphs. 
\subsection{Pseudo-random number generators}
The following text summarised from a patent application publication (see \cite{patent}) clearly explains the application of Ramanujan graphs in the design of pseudorandom number generators.

Pseudorandom numbers may be generated using expander graphs. Expander graphs are a collection of vertices that are interconnected via edges. Generally, a walk around an expander graph is determined responsive to an input seed, and a pseudorandom number is produced based on vertex names. 
Generally, any type of expander graph may be employed to generate pseudorandom numbers. Expander graphs are usually characterized as having a property that enables them to grow quickly from a given vertex to its neighbors and onward to other vertices. An example of a family of graphs that are considered to have good expansion properties are the so-called Ramanujan graphs. 
Expander graphs that are so-called $k$-regular graphs are particularly amenable for use in generating pseudorandom numbers. These $k$-regular graphs are graphs that have the same number of edges emanating from each vertex. The $k$-regular Ramanujan graphs are particularly amenable for use in generating pseudorandom numbers. 
\subsection{Construction of hash functions}
Ramanujan graphs have also been used in constructing secure cryptographic hash functions.  If finding cycles in such graphs is hard, then these hash functions are provably collision resistant.  As an example, there is a specific family of optimal expander graphs for provable collision resistant hash function constructions, namely,  the family of Ramanujan graphs constructed by Pizer (see \cite{Microsoft}).
\begin{figure}
\begin{center}
\includegraphics[width=8cm]{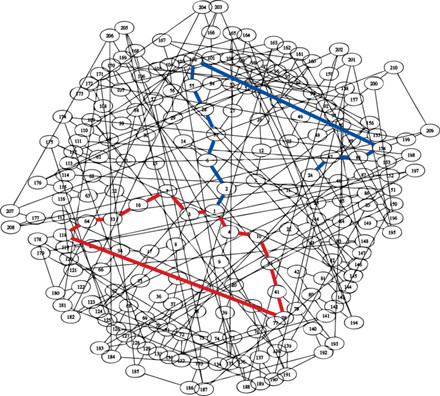}
\caption{Image of a Ramanujan graph that could be used for constructing a cryptograph hash function (for details see \cite{RamHash})}
\end{center}
\end{figure}
\section{Ramanujan machines}
The Ramanujan machine is not actually a machine; it is only a concept. It exists as a network of computers running algorithms and the algorithms outputs conjectures! In his published papers and in his famous Note Books, Ramanujan has sated a large number of results and formulas in the form of conjectures. Most of these conjectures have been proved true by later mathematicians. Because of his propensity to produce conjectures as if from thin air, Ramanujan has been sometimes called the ``conjecture machine'' (see \cite{ConjectureMan}). So it is only appropriate that a machine, whether actual or conceptual, that produces conjectures be called the ``Ramanujan machine''.

Presently what the algorithms do is to come up with probable infinite continued fraction expansions of the constants $e$ and $\pi$. The creators of the machine claim that the algorithm have already generated several such expansions which are conjectured to be true. So, we begin our discussion by explaining what a continued fraction is and then having a look at some of the continued fractions conjectured by Ramanujan. Finally we shall state some of the continued fractions conjectured by the Ramanujan machine.

\section{Continued fractions}
\subsection{Definitions}
An expression of the following form is called a {\em continued fraction}:
$$
a_0+\cfrac{b_1}{ a_1+\cfrac{b_2}{ a_2+\cfrac{b_3}{ a_3+\cdots}}}
$$
It is often expressed in the following form:
$$
a_0 + \cfrac{b_1}{a_1} \pb \cfrac{b_2}{a_2} \pb \cfrac{b_3}{a_3} \pb \cdots
$$
Here $a_0,a_1,a_2,a_3,\cdots$ and $b_1,b_2,b_3,\cdots$ may in general be integers, real numbers, complex numbers, or functions. 

If $a_1,a_2,a_3,\cdots$ are positive integers and $b_1=b_2=b_3=\ldots = 1$ and  then it is called a {\em simple continued fraction}. In a simple continued fraction $a_0$ may be any integer negative, zero or positive. 

A continued fraction may have a finite or an infinite number of terms. For example, we have 
$$
\frac{100}{2.54} = 39 + \cfrac{1}{2} \pb \cfrac{1}{1} \pb \cfrac{1}{2} \pb \cfrac{1}{2} \pb \cfrac{1}{1} \pb \cfrac{1}{4}.
$$
Also, it can be shown that the golden ration $\phi$ has the following infinite  simple continued fraction expansion:
\begin{align*}
\phi 
&= \frac{\sqrt{5}+1}{2}\\
&= 1+ \frac{1}{1}\pb\frac{1}{1}\pb\frac{1}{1}\pb\cdots
\end{align*}
\subsection{Continued fraction expressions for $\pi$ and $e$}
Several continued fraction expressions for $\pi$ are well known. For example, we have the following simple continued fraction expression for $\pi$:
$$
\pi= 3 + \frac{1}{7}\pb \frac{1}{15}\pb \frac{1}{1}\pb\frac{1}{292}\pb\frac{1}{1}\pb\frac{1}{1}\pb\frac{1}{1}\pb\frac{1}{2}\pb\frac{1}{1}\pb\frac{1}{3}\pb\cdots
$$
Unfortunately, the numbers $7,15,1,292,1,\ldots$ have no obvious pattern! But there are continued fraction expressions which are not simple in which the numbers have simple patterns. The oldest such expression is the following one due to William Brouncker (ca. 1660's):
$$
{\pi} = \frac{4}{1}\pb \frac{1^2}{2}\pb \frac{3^2}{2}\pb \frac{5^2}{2}\pb \frac{7^2}{2}\pb \cdots
$$
There is an infinite continued fraction expression for $e$ discovered by Euler (see \cite{Euler, Olds1970}), namely
$$
e=2 + \frac{1}{1}\pb \frac{1}{2}\pb \frac{1}{1}\pb\frac{1}{1}\pb \frac{1}{4}\pb\frac{1}{1}\pb\frac{1}{1}\pb\frac{1}{6}\pb\cdots
$$
\section{Ramanujan and continued fractions}
Ramanujan in his life time published only one continued fraction, the famous Rogers-Ramanujan continued fraction. But his notebooks and the letters he sent to G. H. Hardy before his departure to England contain enough materiel indicating Ramanujan's early interest in continued fractions. These provide ample testimony to the fact that he had studied and used continued fractions extensively (see \cite{RCF1}). 

\subsection{Some early results}
In his first letter to Hardy dated 16 January 1913 Ramanujan has stated the following continued fractions without proof (see \cite{Bruce}). These also appear in his Note Book 2. 
\begin{enumerate}
\item
We have
$$
\displaystyle
\left\{ \frac{\Gamma\left(\frac{x+1}{4}\right)}{\Gamma\left(\frac{x+3}{4}\right)}\right\}^2
=
\frac{4}{x} \pb \frac{1^2}{2x}\pb \frac{3^2}{2x}\pb \frac{5^2}{2x}\pb\frac{7^2}{2x}\pb \cdots
$$
\item
If $ u = x+m+n+1$ and 
$$
P=\frac%
{\Gamma\{\frac{1}{4}u\}\Gamma\{\frac{1}{4}(u-2n)\}
 \Gamma\{\frac{1}{4}(u-2(m-1))\}\Gamma\{\frac{1}{4}(u-2(m+n-1))\}}%
{\Gamma\{\frac{1}{4}(u+2)\}\Gamma\{\frac{1}{4}(u-2m)\}
 \Gamma\{\frac{1}{4}(u-2(n-1))\}\Gamma\{\frac{1}{4}(u-2(m+n))\}}%
$$
then
$$
\frac{1-P}{1+P}=\frac{m}{x}\pb\frac{1^2-n^2}{x}\pb \frac{2^2-m^2}{x}\pb\frac{3^2-n^2}{x}\pb \frac{4^2-m^2}{x}\pb \cdots
$$
\item
If 
\begin{align*}
z&=1+\left(\frac{1}{2}\right)^2x +\left(\frac{1\cdot3}{2\cdot4}\right)^2\cdots\\
y&=\frac{\pi}{2}%
\frac%
{1+\left(\frac{1}{2}\right)^2(1-x) +\left(\frac{1\cdot3}{2\cdot4}\right)^2 (1-x)^2 \cdots}%
{1+\left(\frac{1}{2}\right)^2x +\left(\frac{1\cdot3}{2\cdot4}\right)^2 x^2\cdots}
\end{align*}
then
\begin{align*}
&\frac{1}{(1+a^2)\cosh y} + \frac{1}{(1+9a^2)\cosh 3y} + \frac{1}{(1+25a^2)\cosh 5y} + \cdots\\
&\qquad = 
\frac{1}{2}\frac{z\sqrt{x}}{1}\pb \frac{(az)^2}{1}\pb\frac{(2az)^2 x}{1}\pb \frac{(3az)^2}{1}\pb \frac{(4az)^2x}{1}\pb \cdots
\end{align*}
where $a$ is any quantity.
\item
If
\begin{align*}
u &= \frac{x}{1}\pb \frac{x^5}{1} \pb \frac{x^{10}}{1} \pb \frac{x^{15}}{1}\pb \cdots\\
v &= \frac{\sqrt[5]{x}}{1} \pb \frac{x}{1} \pb \frac{x^2}{1} \pb \frac{x^3}{1}\pb \cdots
\end{align*}
then
$$
v^5 =u\cdot\frac{1-2u+4u^2-3u^3+u^4}{1+3u+4u^2+2u^3+u^4}.
$$
\item
We have 
$$
\displaystyle
\frac{1}{1} \pb \frac{e^{-2\pi}}{1} \pb \frac{e^{-4\pi}}{1} \pb \frac{e^{-6\pi}}{1}\pb \cdots 
=
\left( \sqrt{\frac{5+\sqrt{5}}{2}} - \frac{\sqrt{5}+1}{2}\right)\sqrt[5]{e^{2\pi}}.
$$
\end{enumerate}

\subsection{Rogers-Ramanujan continued fraction}

The Rogers-Ramanujan continued fraction was first discovered by Rogers and published with proof in 1888 and then rediscovered independently and published without proof by Ramanujan in 1914. If we write
\begin{align*}
G(q)
& = \sum_{n=0}^\infty \frac{q^{n^2}}{(1-q)(1-q^2)\cdots (1-q^n)}\\
H(q)
& = \sum_{n=0}^\infty \frac{q^{n^2+n}}{(1-q)(1-q^2)\cdots (1-q^n)}
\end{align*}
then the following continued faction is known as the Rogers-Ramanujan continued fraction: can 
\begin{align*}
R(q) 
& = q^{1/5}\frac{H(q)}{G(q)}\\
&= \frac{q^{1/5}}{1}\pb\frac{q}{1}\pb\frac{q^2}{1}\pb\frac{q^3}{1}\pb\cdots
\end{align*}
In his very first letter to Hardy, Ramanujan has stated several properties of the function $R(q)$ (see \cite{RCF2}).
\section{Conjectures by Ramanujan machine}
The Ramanujan machine was initially designed to output conjectures  regarding continued fraction expressions for the mathematical constants $\pi$ and $e$. The machine has produced several such conjectures and most of the early ones have been proved to be correct. The capability of the machine has now been expanded to generate such expressions for other constants like the Catatlan constant, values of the Reimann zeta function, etc. 

``Fundamental constants like $e$ and $\pi$ are ubiquitous in diverse fields of science, including physics, biology, chemistry, geometry, and abstract mathematics. Nevertheless, for centuries new mathematical formulas relating fundamental constants are scarce and are usually discovered sporadically by mathematical intuition or ingenuity.'' (see \cite{RamanujanMachineWebsite})

\begin{enumerate}
\item
Conjectures for $\pi$

The following is a conjecture involving $\pi$ discovered by the ramanujan machine and which has been proved to be true.

$$
\frac{4}{3\pi - 8}=3 +\frac{-1\cdot 1}{6}\pb \frac{-2\cdot 3}{9}\pb \frac{-3\cdot 5}{12}\pb\frac{-4\cdot 7}{15}\cdots
$$

\item
Conjectures for $e$

Here is an elegant continued fraction representation for $e$  discovered by the Ramanujan machine.
$$
e=3 + \frac{-1}{4}\pb\frac{-2}{5}\pb \frac{-3}{6}\pb \frac{-4}{7}\pb\cdots
$$

\item
Conjectures for $\log 2$:

The following conjecture involving $\log 2$ was discovered after 1 May 2020 and has not yet been proven.
\begin{align*}
\frac{1}{1-\log 2}
& = a_0 + \frac{a_1}{b1} \pb \frac{a_2}{b_2}\pb \cdots \pb \frac{a_n}{b_n}\pb\cdots\\
& = 4 + \frac{-8}{14}\pb \frac{-72}{30}\pb\cdots
\end{align*}
where
$$
a_n=3n^2+7n+4, \quad b_n=-2n^2(n+1)^2.
$$
\item
Conjectures on Catalan's constant

The Catalan's constant $G$ is defined by
$$
G=\frac{1}{1^2}-\frac{1}{3^2}+\frac{1}{5^2}-\frac{1}{7^2}=\cdots
$$
It is not known whether G is irrational, let alone transcendental. There is Ramanujan connection to this constant. Ramanujan has extensively studied the following inverse tangent integral function
$$
\text{Ti}(x)=\int_0^x\frac{\tan^{-1} x}{x}\, dx.
$$
It can be shown that $G=\text{Ti}(1)$.
The Ramanujan machine has produced several conjectures regarding $G$ most of whichh have not yet been proven. The following is a typical conjecture on $G$ discovered by the  machine:
\begin{align*}
\frac{1}{2G}
& = a_0 + \frac{a_1}{b_1} \pb \frac{a_2}{b_2}\pb \cdots \pb \frac{a_n}{b_n}\pb\cdots\\
& = 1 + \frac{-2}{7}\pb \frac{-32}{19}\pb\cdots
\end{align*}
where
$$
a_n=-2n^4\quad b_n=3n^2+3n+1.
$$
\item
Conjectures on Apery's constant

The Reimann zeta function, $\zeta(s)$, is defined by 
$$ 
\zeta(s)=\frac{1}{1^s}+\frac{1}{2^s}+\frac{1}{3^s}+\cdots
$$
The number $\zeta(3)$ is defined as the Apery's constant. It is named after Roger Apéry (1916 - 1994), a Greek-French mathematician. This constant appears in several physical problems  including in the study of quantum electrodynamics.

The Ramanujan machine has generated several conjectures regarding the Apery's constant. For example, 
$$
\frac{1}{\zeta(3)}= a_0 + \frac{a_1}{b_1} \pb \frac{a_2}{b_2}\pb \cdots \pb \frac{a_n}{b_n}\pb\cdots
$$
where
$$
a_n = -n^6, \quad b_n = n^3+(n+1)^3.
$$
\end{enumerate}
\section{Ramanujan in digital signal processing}
First some background regarding digital signal processing. A signal is an electrical or electromagnetic current that is used for carrying data from one device to another. 
It is the key component behind virtually all communication computing, networking and electronic devices.
Any quality, such as a physical quantity that exhibits variation in space or time can be used as a signal. A signal can be audio, video, speech, image, sonar and radar-related and so on. 

Signal processing is a subfield of electrical engineering that deals with analysing, modifying, and synthesizing signals such as sound, images, and scientific measurements. Digital signal processing (DSP) is the process of analyzing and modifying a signal to optimize or improve its efficiency or performance. It involves applying various mathematical and computational algorithms to analog and digital signals to produce a signal that is of higher quality than the original signal. Digital signal processing is primarily used to detect errors, and to filter and compress analog signals in transit.
\section{Ramanujan sum}
\subsection{Definiton}
In a paper published in 1918 (see \cite{R1918}), Ramanujan introduced the following sum where $n,d,k$ are positive integers and where  $(k,d)$ denotes the greatest common divisor of $k$ and $d$:
$$
c_q(n) = \sum_{k=1\atop (k,q)=1}^n \cos \frac{2kn\pi}{q}.
$$
The expression on the right can be rewritten in the following form:
$$
c_q(n)=\sum_{k=1\atop (k,q)=1}^n e ^{\frac{2kn\pi}{q}i}.
$$
The complex number $e^{-\frac{2\pi}{q}i}$ is a $q$-th root of unity and it is sometimes denoted by $W_q$. Using this notation, we also have
$$
c_q(n)= \sum_{k=1\atop (k,q)=1}^n W_q^{-kn}.
$$
The sum $c_q(n)$ is now generally referred to as the {\em Ramanujan sum}.

For example let $q=6$. The integers $k$ lying in the range $1,2,3,4,5,6$ and satisfying $(k,q)=1$ are $1,5$ only. Hence we have
\begin{align*}
c_6(n)& = \cos \frac{4\pi n}{6} + \cos \frac{10\pi n}{6}\\
&=e^{\frac{4\pi n}{6}i}+e^{\frac{10\pi n}{6}i}.
\end{align*}
A few values of $c_6(n)$ are given in Table \ref{c6n}.  It may be noted that $c_6(n+6)=c_6(n)$ for all $n$.
\begin{table}[!h]
\begin{center}
\begin{tabular}{c|ccccccccccccc}
\hline
$n$ & 0 & 1 & 2 & 3 & 4 & 5 & 6 & 7 & 8 & 9 & 10 & 11 & $\cdots$ \\
\hline
$c_6(n)$ & $2$ & $1$ & $-1$ & $-2$ & $-1$ & $1$ & $2$ & $1$ & $-1$ & $-2$ & $-1$ & $1$ & $\cdots$\\
\hline
\end{tabular}
\caption{Values of $c_6(n)$}\label{c6n}
\end{center}
\end{table}
\subsection{Some properties}
A few of the many interesting properties of the Ramanujan sum are listed below.
\begin{enumerate}
\item
The Ramanujan sum $c_q(n)$ is a period function of $n$ with period $q$, that is,
$$
c_q(n+q)=c_q(n)\quad \text{for all $n$}.
$$
\item
$c_q(n)$ is an integer valued function.
\item
The function $c_q(n)$ is a multiplicative function of $q$, that is, if $(q_1,q_2)=1$ (in other words, if $q_1$ and $q_2$ are relatively prime), then
$$
c_{q_1q_2}(n)=c_{q_1}(n)c_{q_2}(n).
$$
\item
Using the Mobius function $\mu(n)$ there is an explicit formula for $c_q(n)$:
$$
c_q(n)=\sum_{d\atop(\text{$d$ divisor of $(q,n)$})} \mu(q/d)d.
$$
Note that the Mobius function is defined as follows:
$$
\mu(n)=
\begin{cases}
1 & \text{if $n=1$}\\
(-1)^k & \text{if $n$ is  a product of $k$ distinct prime numbers}\\
0&\text{otherwise}
\end{cases}
$$
\item
Let $q_1\ne q_2$ and let $l$ be any common multiple of $q_1$ and $q_2$. Then
$$
\sum_{n=0}^{l-1}c_{q_1}(n)c_{q_2}(n)=0.
$$
\item
The Euler's totient function $\phi(n)$ is an arithmetic function which takes the number of positive integers not exceeding $n$ and relatively prime to $n$ as its value. Then
\begin{itemize}
\item
$
\displaystyle{\lim_{x \rightarrow \infty}\frac{1}{x}}\sum_{n\le x}c_r(n)c_s(n)
=
\begin{cases}
c_r(h)&\text{if $r=s$}\\
0 &\text{otherwise}
\end{cases}
$
\item
$
\displaystyle{\lim_{x \rightarrow \infty}\frac{1}{x}}\sum_{n\le x}c_r(n)c_s(n+h)
=
\begin{cases}
\phi(r)&\text{if $r=s$}\\
0 &\text{otherwise}
\end{cases}
$
\end{itemize}
\end{enumerate}
\subsection{Ramanujan-Fourier series}
An arithmetic function is a function defined on the set of positive integers which takes real or complex numbers as values. For example, the {\em divisor function} $d(n)$ whose value is the number of divisors of $n$ is an arithmetic function.  Ramanujan's principal objective in studying the properties of the function $c_q(n)$ was to obtain an expression for an arithmetic function $x(n)$    in the form of a series as follows.
$$
x(n) = \sum_{q=1}^\infty a_q c_q(n).
$$
A series of this form is called {\em Ramanujan series} (or {\em Ramanujan expansion} or some-
times {\em Ramanujan-Fourier series} or {\em Ramanujan-Fourier Transform}) (see \cite{Murty}).

As an example, Ramanujan discovered the following expansion for the divisor function $d(n)$:
$$
d(n) = - \sum_{q=1}^\infty \frac{\log q}{q} c_q(n).
$$
As another example, consider the arithmetic function $\sigma(n)$ which takes the sum of the divisors of $n$ as its value. Thus $$
\sigma(6)=1+2+3+6=12, \quad \sigma(7)=1+7=8.
$$
Ramanujan obtained the following series expansion for the function  $\sigma(n)$:
$$
\sigma(n)=\frac{\pi^2n}{6}\sum_{q=1}^\infty \frac{c_q(n)}{q^2}.
$$
\subsection{An early application of Ramanujan sums to signal processing}
The discrete and fast Fourier transforms are suited for the analysis of periodic and quasiperiodic sequences. But they are not suitable to discover the features of aperiodic
sequences such as low-frequency noise. The Ramanujan-Fourier transform has been applied to examine whether there are any hidden patterns in such sequences. The general idea is to compare the observed patterns and properties in experimental data with the properties of the Ramanujan-Fourier transforms of well known arithmetical functions in number theory and to see whether there is any match between the two sets of properties. 

This idea was tested using certain observed data from the study of black holes in astronomy and the general 
conclusion was that many of these
processes may be described using prime number theory. It has also been suggested that Ramanujan-Fourier transform could also be used in the study of  radio-frequency oscillators close to phase locking (see \cite{Planat}).

\section{Ramanujan spaces}
For a fixed positive integer $q$, consider the  matrix:
$$
B_q
=
\begin{bmatrix}
c_q(0) & c_q(q-1) & c_q(q-2) & \cdots & c_q(1)\\
c_q(1) & c_q(0) & c_q(q-1) & \cdots & c_q(2)\\
c_q(2) & c_q(1) & c_q(0) & \cdots & c_q(3)\\
\vdots &        &        &        &   \\
c_q(q-1) & c_q(q-2) & c_q(q-3) & \cdots & c_q(0)
\end{bmatrix}
$$
Let $\mathbb C^q$ be the $q$-dimensional complex vector space whose vectors are represented as column vectors. The subspace of $\mathbb C^q$ generated by the columns of $B_q$ is called the {\em Ramanujan space} (or sometimes, the {\em Ramanujan subspace}) and is denoted by $S_q$. It is shown that the space $S_q$ has dimension $\phi(q)$ where $\phi$ is the Euler's  totient function. Also, any $\phi(q)$ consecutive columns in $B_q$, in particular the first $\phi(q)$ columns of $B_q$, form a basis for $S_q$. 

As a concrete example, let us consider $B_6$ using the values of $c_6(n)$ given in Table \ref{c6n}.
$$
B_6
=
\begin{bmatrix}
\phantom{-}2 & \phantom{-}1 & -1 & -2 & -1 & \phantom{-}1 \\
\phantom{-}1 & \phantom{-}2 & \phantom{-}1 & -1 & -2 & -1 \\
-1 & \phantom{-}1 & \phantom{-}2 & \phantom{-}1 & -1 & -2 \\
-2 & -1 & \phantom{-}1 & \phantom{-}2 & \phantom{-}1 & -1 \\
-1 & -2 & -1 & \phantom{-}1 & \phantom{-}2 & \phantom{-}1 \\
\phantom{-}1 & -1 & -2 & -1 & \phantom{-}1 & \phantom{-}2 
\end{bmatrix}
$$
It may be noted that $\phi(6)=2$ and  the first two columns are linearly independent, and further any other column is a linear combination of the first two columns. 

We defined $S_q$ as a $q$-dimensional subspace of vectors in
$\mathbb C^q$, spanned by
columns of $B_q$. However, the
elements of $S_q$ are sometimes regarded as period-$q$ sequences. So, when
we say that $x(n)\in S_q$, it is understood that the vector of size $q$
is extended periodically to obtain the sequence $x(n)$.

The Ramanujan spaces has an important periodicity property. A discrete-time signal $x(n)$ is said to be periodic if there exists an integer $R$ such that
$$
x(n) = x(n + R)
$$
for all $n$, and the integer $R$ is called a repetition interval. The period $P$ (an integer) is the smallest
positive repetition interval. It can be shown that any repetition interval is an integer multiple
of $P$. If $x_1 (n)$ and $x_2 (n)$ have periods $P_1$ and $P_2$, their
sum $x_1 (n) + x_2 (n)$ has a repetition interval $R = \text{lcm}\, (P_1, P_2 )$, so its period $P$ is either this lcm or a
proper divisor of it.

It has been shown that all vectors in the Ramanujan space $S_q$ has period exactly equal to $q$, in particular, it cannot be smaller than $q$. More generally, let us consider a set of $K$ Ramanujan spaces $S_{q_m}$, $m=1,\ldots,K$ and $x_m(n)\in S_{q_m}$. Let $N$ be the lcm of $q_1, \ldots, q_K$. Then the vector
$$
x(n)=\sum_{m=1}^K x_m(n)
$$
is periodic with period exactly equal to $N$; in particular, it cannot be smaller (see \cite{Vaidya})

Let $x(n)$ be finite sequence of length. Let us assume that the sequence is extended for all values of $n$ by stipulating that
$$ 
x(n+N)=x(n), \text{ for all $n$.}
$$
It can be shown that such a sequence, considered as a vector $\mathbf x$ in $\mathbb C^N$, can be expressed as a sum of vectors in the Ramanujan spaces $S_{q_1}, S_{q_2}, \ldots, S_{q_K}$ where $q_1, q_2, \ldots, q_K$ are the divisors of $N$:
$$
\mathbf x = \sum_{q_i|N} \mathbf x_{q_i}, \text{ where $\mathbf x_{q_i} \in S_{q_i}$}.
$$
This has been called the {\em Ramanujan FIR Representation}\footnote{FIR is an acronym for Finite Impulse Response.} of the signal $x(n)$. These representations can be used to determine hidden periodicities in signals and also for denoising of signals. For more details the reader may refer to \cite{Vaidya2}. 

\subsection{Applications}
The following quote from the concluding section of Vaidyanathan's paper on Srinivasa Ramanujan and
signal-processing problems (see \cite{Vaidya3}) gives an indication of the various fields where the ideas of Ramanujan sums and Ramanujan subspaces have been successfully applied.

''In this paper, we presented an overview of the impact of Ramanujan sums in signal processing,
especially integer-period estimation in real or complex signals. These methods have recently been
used in identification of integer periodicities in DNA molecules and in protein molecules. These new methods are quite competitive and often work better than
existing state of the art methods. Ramanujan filter banks have also been shown to be applicable in
the identification of epileptic seizures in patients, which are characterized by sudden appearance
of periodic waveforms in the measured EEG records. Other interesting applications have
recently been reported by a number of authors such as, for example in source-separation,
RF communications, ECG signal processing and brain-computer interfacing. More
recently, a well-known algorithm called the MUSIC algorithm, which is popularly used for
identifying sinusoids in noise, has been extended to the case of integer period identification using
Ramanujan-subspace ideas. This method, known as iMUSIC, has also been compared with
other well-known methods for multipitch estimation.'' 
\section{Conclusion}

``It is satisfying indeed when one finds that a well-known mathematical concept has practical
impact in engineering, even though the original mathematical ideas may not have been inspired
by any such application. Engineers have seen this happening over and over again in disciplines
such as information theory, coding, digital communications, system theory and machine learning.
Indeed, the classical view that pure mathematics of the highest quality is bound to be `useless'
for real-life applications is evidently not valid as the last several decades of science and
engineering have amply demonstrated. What is often regarded as pure mathematics sometimes
impacts engineering in wonderful ways. A classic example is the theory of finite fields and
rings, which has impacted the practice of error-correction coding in digital communications, data
compression, and digital storage. Other examples include graph theory which has had many
engineering, network and signal processing applications. Yet another is number theory which has
impacted nearly all aspects of science and engineering such as acoustical hall designs, computer
music, and so forth.'' (see \cite{Vaidya3})

%
%

%
\end{document}